\documentclass[12pt]{article}
\usepackage{graphicx}
\usepackage{amsmath}
\usepackage{graphics}
\usepackage{epsfig}
\usepackage{amssymb}

\usepackage[authoryear]{natbib}
\newcommand{\be}{\begin{eqnarray}}
\newcommand{\ee}{\end{eqnarray}}
\newcommand{\bea}{\begin{eqnarray*}}
\newcommand{\eea}{\end{eqnarray*}}

\newcommand{\BLUE}{\scalebox{0.5}{$\mathrm{BLUE}$}}
\newcommand{\SLSE}{\scalebox{0.5}{$\mathrm{SLSE}$}}
\newcommand{\LSE}{\scalebox{0.5}{$\mathrm{LSE}$}}

\newcommand{\ds}{\displaystyle}

\newtheorem{theorem}{Theorem}[section]
\newtheorem{lemma}{Lemma}[section]
\newtheorem{proposition}{Proposition}[section]
\newtheorem{corollary}{Corollary}[section]

\newtheorem{remark}{Remark}[section]


\begin{document}

\title{Optimal designs for regression models with autoregressive errors structure}

\author{Holger Dette  \\
Fakult\"at f\"ur Mathematik \\
Ruhr-Universit\"at Bochum \\
44780 Bochum \\
Germany
\and
Andrey Pepelyshev, Anatoly Zhigljavsky \\
School of Mathematics\\
 Cardiff  University\\
  Cardiff, CF24 4AG\\
  UK
}
\date{}
\maketitle

\begin{abstract}

In the one-parameter regression model  with AR$(1)$ and AR$(2)$  errors
we find   explicit expressions and a continuous approximation of the optimal discrete design for the
signed least square estimator. The results are used to derive the optimal variance of the best linear
estimator in the continuous time model and  to  construct efficient estimators and
corresponding optimal designs for finite samples. The resulting procedure (estimator and design)
provides  nearly the same efficiency as the weighted least squares and its variance
is close to the optimal variance   in the
continuous time
model. The results are illustrated by several examples demonstrating the feasibility of our approach.
\end{abstract}

Keywords and Phrases: linear regression; correlated observations; signed measures; optimal design; BLUE; AR  processes; continuous autoregressive model

AMS Subject classification: Primary 62K05; Secondary 31A10

% \tableofcontents

\section{Introduction}
\label{sec1}
\def\theequation{1.\arabic{equation}}
\setcounter{equation}{0}

%\subsection{Models with errors satisfying discrete AR$(p)$ process}

Consider a linear   regression model
\be
 y_j=\theta^Tf(t_j)+\epsilon_j\, \;\; (j=1, \ldots, N)\, ,
 \label{eq:model}
\ee
where $\theta \in \mathbb{R}^m$ is a vector of unknown parameters, $f(t)=(f_1(t), \ldots, f_m(t))^T$ is a vector  of linearly independent functions defined
on some interval, say $[A,B]$, and $\epsilon_1, \ldots , \epsilon_N$ are  random errors with $\mathbb{E}[\epsilon_j]=0$ for all $j=1, \ldots, N$ and
covariances  $\mathbb{E}[\epsilon_j \epsilon_k] =\rho(t_j-t_k)$.  It is well known
that  the use of optimal or efficient designs yields to a reduction of costs  by a statistical inference  with a minimal number of experiments without loosing any accuracy.
Optimal design theory has been studied intensively for the case when errors are uncorrelated using tools from convex optimization theory  [see  \citep{pukelsheim2006}], but
the design problem in the case of  dependent data is substantially harder because the  corresponding optimization problems are usually non-convex.
Most authors use asymptotic arguments to construct optimal designs, which do not solve the problem of non-convexity
 [see for example  \cite{sackylv1966,sackylv1968}, \cite{bickherz1979}, \cite{N1985a}, \cite{zhidetpep2010}, \cite{dette2015design}].
Some optimal designs for the location  model (in this case the optimization problems are in fact convex) and  for a few one-parameter
linear models have been discussed in  \cite{bolnat1982}, \cite{N1985a}, Ch.\ 4, \cite{naether1985b},  \cite{pazmue2001}   and  \cite{muepaz2003} among others].
{Recently, for multi-parameter models, \cite{DetPZ2012} determined a necessary condition for the optimality of (asymptotic) designs
for least squares estimation.
\cite{DPZ2014universal} studied nearly universally optimal designs,
while \cite{DetPZ2016} constructed new matrix-weighted estimators with corresponding optimal designs,
which are very close to the best linear unbiased estimator with corresponding  optimal designs.}
Although these results are promising, they rely on certain
structural assumptions on the covariance kernel. For example,  \cite{DetPZ2012} assume
that the regression functions in model \eqref{eq:model}  are eigenfunctions of an integral operator associated with the covariance kernel of the error process
and  \cite{DetPZ2016}  assume that the covariance kernel is triangular [see \cite{mehr1965certain} for an exact definition].  While these results cover
the frequently used AR($1$)-process as error structure, they are not applicable in models with autoregressive error processes of
  larger order.

The goal of the present paper is to give first insights in the optimal design problem for linear regression models  with autoregressive error processes. We concentrate on a one-parameter linear regression model with an AR$(1)$ and AR$(2)$-error process. In Section \ref{sec2} we will introduce a signed least squares estimator and consider approximate designs on the design space ${\cal T}=\{t_1,\ldots,t_N\}$, where the weights are not necessarily non-negative.
We  determine the optimal (signed) approximate design for signed least squares estimation, such that the signed least squares estimator has the same variance as the weighted least squares estimator based on observations at the experimental conditions $t_1,\ldots,t_N$. In Section \ref{sec3} we consider the one-parameter linear regression model with
autoregressive errors of order $1$ and study the asymptotic behavior of the signed least squares estimator with corresponding optimal design as the sample size tends to infinity. Section \ref{sec4} is devoted to the case of an AR$(2)$-error process, where the situation is substantially more complicated. These results are then   used in Section \ref{sec5}, where we consider the problem of constructing designs for signed least squares estimation in finite sample situations. We provide a procedure such that the signed least squares estimator with corresponding ``optimal'' design has nearly the same efficiency as the weighted least squares estimator with corresponding optimal design. Finally, the results are illustrated in several numerical examples.

\section{Various least squares estimators }
\label{sec2}
\def\theequation{2.\arabic{equation}}
\setcounter{equation}{0}

 For estimating $\theta$, we  use the following two estimators:  the best linear unbiased estimator (BLUE)
\bea
%\label{blue}
 \hat \theta_{\BLUE,N}=(X^T\Sigma^{-1}X)^{-1}X^T\Sigma^{-1}Y
\eea
and the signed least squares estimator (SLSE)
\be
 \widehat \theta_{\SLSE,N}  = ({X}^T\textbf{S}{X})^{-1}{X}^T\textbf{S} {Y},
 \label{eq:slse}
\ee
where $X=(f_i(x_j))_{j,i=1}^{N,m}$ is the design matrix of size $N \times m$,
$\textbf{S}$ is an $N\!\times\! N$ diagonal matrix with entries $+1$ and $-1$ on the diagonal and
$\Sigma=(\rho(t_i-t_j))_{i,j=1}^N$ is the covariance matrix of observations. If $\textbf{S}$ is the $N\!\times\! N$ identity matrix, then
SLSE coincides with the ordinary least squares estimator (LSE).
The covariance matrix of the BLUE and the  SLSE are given by
\bea
\label{blue_var}
 \mathrm{Var}(\hat \theta_{\BLUE,N})&=&(X^T\Sigma^{-1}X)^{-1}~, \\
 \mathrm{Var}(\hat \theta_{\SLSE,N}) &=&({X}^T\textbf{S}{X})^{-1} (X^T\textbf{S}\Sigma^{}\textbf{S}X) ({X}^T\textbf{S}{X})^{-1},
\eea
respectively. Throughout this paper  we concentrate on the  one-parameter regression model
\be \label{one}
y_j=\theta f(t_j)+\epsilon_j,
\label{eq:model-onepar}
\ee
and remark that an extension to the multi-parameter model \eqref{eq:model} could be performed following the discussion
in \cite{DetPZ2016}. A  design  on the (fixed) design space ${\cal T}=\{t_1,\ldots,t_N\}$
is an arbitrary discrete signed measure of the form
$ \xi = \{t_1,\ldots,t_N;w_1,\ldots,w_N\}$, where $w_i=s_i p_i$, $s_i \in \{-1,1\}$, $p_i \ge 0$, $i=1, \dots, N$,
and \mbox{$\sum_{i=1}^N p_i=1$}.  The variance of the SLSE for the design $\xi$ is given by
\be
\label{eq:D-SLSw}
D (\xi)= \mathrm{Var}(\hat \theta_{\SLSE,N})=\sum^N_{i=1}\sum^N_{j=1} \rho(t_i-t_j) w_i w_j f_i f_j\Big/\Big(\sum^N_{i=1} w_i f_i^2\Big)^2  \, ,
\ee
where we use the notation $f_i=f(t_i)$  throughout this paper.
The optimal design problem  consists in the minimization
of this expression with respect to the weights $w_1,\dots,w_N$ assuming that the observation points $t_1,\dots,t_N$ are fixed.
 Despite the fact that the functional $D$ in \eqref{eq:D-SLSw} is not  convex as a function of  $w_1,\dots,w_N$,
the problem of determining the optimal weights can be easily solved by a simple application of the  Cauchy-Schwarz inequality.
The proof of the following lemma is given in \cite{DetPZ2016}; see also  Theorem 5.3 in  \cite{N1985a}, where
 this result was proved in  a  slightly different form.

\begin{lemma}
\label{lem:2-1}
Assume that the matrix ${\Sigma}=(\rho(t_i-t_j))_{i,j=1,\ldots,N}$
is positive definite and \mbox{$f_i\neq 0$} for all $i=1,\ldots,N$. Then
the optimal weights $w_1^*,\dots,w_N^*$ minimizing the expression \eqref{eq:D-SLSw} are given by
\be
\label{eq:weights}
 w_i^*= \,{\mathbf{e}_i^T {\Sigma}^{-1}\mathbf{f}}/{f_i}; \qquad i=1,\dots,N,
\ee
where $\mathbf{f}=(f_1,\ldots,f_N)^T$, $\mathbf{e}_i=(0,\ldots,0,1,0, \ldots,0)^T \in \mathbb{R}^N$ is the $i$-th unit vector.
Moreover, for the design $\xi^*=\{t_1,\ldots,t_N;w_1^*,\ldots,w_N^*\}$ with weights \eqref{eq:weights} we have
$D (\xi^*)=D^*$, where $D^*=1/(\mathbf{f}^T\,{\Sigma}^{-1}\mathbf{f})$ is the variance of the BLUE defined in  \eqref{blue_var}.
\end{lemma}
Note that the optimal weights in Lemma  \ref{lem:2-1} are not uniquely defined. In fact, they can always be multiplied by a constant without changing their
optimality. In the following discussion we will consider the case where the points $t_i$  are given by the equidistant points  on the interval $[A,B]$
and the sample size $N$ tends to infinity. Heuristically the BLUE converges in this case to the BLUE in the continuous time model,
where the full trajectory of the stochastic process can be observed. 
Note that for any finite  $N$ the SLSE with the optimal weights defined in Lemma \ref{lem:2-1} has the same variance as the BLUE.

Further we study the asymptotic properties  of the SLSE and the optimal weights $w^*_i$ defined in \eqref{eq:weights} as the sample size increases.
 In many cases we will be able to
approximate an $N$-point design  $ \xi = \{t_1,\ldots,t_N;w_1^*,\ldots,w_N^*\}$ with optimal weights defined
 in \eqref{eq:weights} by a signed measure (an approximate design) of the form
\be
\label{cont_design}
 \xi(dt)=P_A\delta_A(dt)+P_B\delta_B(dt)+p(t)dt\, ,
% \xi(dt)=Q_A\frac{\delta_A(dt)-\delta_{A+\epsilon}(dt)}{\epsilon}
% +Q_B\frac{\delta_B(dt)-\delta_{B-\epsilon}(dt)}{\epsilon}+ P_A\delta_A(dt)+P_B\delta_B(dt)+p(t)dt\, ,
\ee
where $\delta_A(dt)$ and $\delta_B(dt)$ are Dirac-measures concentrated at  the point $A$ and $B$, respectively, and
$p(\cdot)$ is a density function (not necessarily non-negative) on the interval  $[A,B]$.
%We will normalize $\xi(dt)$ so that the total mass is 1; that is, $|\xi|([A,B])= |P_A|+|P_B|+\int_A^B |p(t)|dt=1$.
Approximate designs  of the from \eqref{cont_design} are easier to understand and analyze than discrete designs
of the form $ \xi = \{t_1,\ldots,t_N;w_1^*,\ldots,w_N^*\}$, and we will illustrate in Section \ref{sec3} and \ref{sec4}
the derivation of the limits in the case of autoregressive error processes of order one and two, respectively.  \\
As already mentioned in the introduction the AR$(1)$ process corresponds to a triangular kernel and could also be
treated with methodology developed in \cite{DetPZ2016}. We discuss it here because for this case the arguments
are simpler than for the AR$(2)$. In fact, for the AR(2) error process the  derivation of asymptotically   optimal weights $w_1^*,\ldots,w_N^*$
of the form \eqref{eq:weights} as the sample size tends to infinity is substantially harder  and
we have to slightly modify the continuous  approximation of the   form  \eqref{cont_design}
[see Section \ref{sec4} for more details].

\section{ Autoregressive errors of order one}
\label{sec3}
\def\theequation{3.\arabic{equation}}
\setcounter{equation}{0}

Consider the regression model \eqref{eq:model} with $N$  equidistant points
\be
\label{points}
t_j =A+(j-1)\Delta~,~ (j=1, \ldots, N)
\ee
on the interval $[A,B]$, where  $\Delta=(B-A)/(N-1)$.
Assume that the errors $\epsilon_1, \ldots, \epsilon_N$ in \eqref{one} satisfy the discrete AR(1) equation
\be
\label{AR1}
\epsilon_j-a\epsilon_{j-1}=z_j
\ee
for some $0<a<1$, where $\epsilon_1\sim N(0,\sigma^2)$ and $z_2, \ldots, z_N$ are Gaussian independent identically distributed random variables
with mean 0 and variance $\sigma_z^2=(1-a^2)\sigma^2$. Without loss of generality, we assume $\sigma^2=1$.

\begin{remark} \label{arparameter} {\rm
Note that discrete AR$(1) $ processes  \eqref{AR1} are usually considered for the parameter  $-1<a<1$.
For the subsequent discussion we need a continuous analogue, say $\{ \varepsilon (t) \}_{t  \in [A,B]}$, of  the discrete AR$(1)$ error process,
which is in fact available in the case $0 <a<1$; see \cite{CT1987}. The corresponding process
with drift is denoted  by $y(t) = \theta f(t) + \varepsilon (t)$,  ${t  \in [A,B]}$.
However, for $-1<a < 0 $  the discrete AR$(1)$ process \eqref{AR1} does not have a continuous real-valued analogue
and  therefore in this case the limiting behavior of our estimators and designs is much harder to understand. \\
It is also worthwhile to mention that the autocovariance function of errors $\epsilon_1, \ldots, \epsilon_N$ is
given by
$$
   \mathbb{E}[\epsilon_j \epsilon_k] =  \rho(t_j-t_k)  = e^{- \lambda | t_j-t_k |} = e^{ \lambda  t_j }e^{- \lambda  t_k}~~
   \mbox{ if }  t_j \leq  t_k  ,
   $$
 where $\lambda=-\ln (a) / \Delta$.  Thus, if $a \in (0,1)$, the AR$(1)$ error process has a triangular covariance kernel
 in the sense of \cite{mehr1965certain}, and the results of \cite{DetPZ2016} are applicable.  In the following discussion
 we provide a different derivation of the asymptotically optimal weights, because the arguments
 will be useful for the discussion of an AR$(2)$-error process in Section \ref{sec4}.
}
\end{remark}

\medskip

For an AR$(1)$-error process, the inverse of the covariance matrix $\Sigma=(\rho(t_i-t_j))_{i,j=1}^N$ is given by the  tridiagonal matrix
\bea
 \Sigma^{-1}=\frac{1}{S}
 \begin{pmatrix}
 1&k_{1}&0&0  &\ldots\\
 k_1 & k_0    &k_1&0&\ldots\\
 0& k_1    &k_0&k_1&0\\
 \vdots&\ddots&\ddots&\ddots&\ddots&\ddots\\
 && 0    &k_1&k_0&k_1\\
 &&0& 0&k_{1}&1\\
 \end{pmatrix}~,
\eea
where  $k_0=1+a^2=1+e^{-2\lambda\Delta}$, $k_1=-a=-e^{-\lambda\Delta}$, $S=1-a^2=1-e^{-2\lambda\Delta}$ and $\lambda=-\ln (a) / \Delta$.
Recalling the definition of the optimal  weights $w_i^* $, $i=2,\ldots, N-1,$  in  \eqref{eq:weights} we have
\bea
 S{w_i^*}{}f(t_i)&=&k_1f_{i-1}+k_0f_i+k_1f_{i+1}=(1+a^2)f_i-af_{i-1}-af_{i+1}\\
 &=&a(2f_i-f_{i-1}-f_{i+1})+(1-2a+a^2)f_i\\
 &=&a(2f_i-f_{i-1}-f_{i+1})+(a-1)^2f_i\, .
\eea
We now assume that $\lambda=-\ln (a) / \Delta$  is fixed and $\Delta  =(B-A)/(N-1) \to 0$. Since $S(\Delta)=S'(0)\Delta+o(\Delta)$ with $S'(0)=2\lambda$ and $a=1-\lambda\Delta+o(\Delta)$, we obtain
\bea
 {w_i^*}{}f(t_i)&=&\frac{\Delta}{S(\Delta)}\cdot\frac{a(2f_i-f_{i-1}-f_{i+1})+(a-1)^2f_i}{\Delta^2}\Delta\\
 &=&\frac{1}{S'(0)}[-f''(t_i)+\lambda^2f(t_i)]\Delta+o(\Delta).
\eea
Thus, we have
\bea
 \frac{w_i^*}{\Delta}=\frac{1}{2\lambda f(t_i)}[-f''(t_i)+\lambda^2f(t_i)]+O(\Delta).
\eea
Therefore, for small $\Delta$, the discrete signed measure $\{t_2,\ldots,t_{N-1} ;w_2^*,\ldots,w_{N-1}^*\}$
is approximated by the continuous signed measure with density
\be
 p(t)=-\frac{1}{2\lambda f(t)}\Big(f''(t)-\lambda^2f(t)\Big).
\label{eq:ar1-pt}
\ee
Now we consider the weights at the boundary points.
For the left boundary weight, we obtain
\bea
 {w_1^*}{}f(t_1)&=&
 \frac{f_1+k_1f_2}{S(\Delta)}=\frac{\Delta}{S(\Delta)}\cdot\frac{f_1-af_2}{\Delta}\\
 &=&\frac{\Delta}{S(\Delta)}\Big[\frac{f_1-f_2}{\Delta}+\frac{f_2-af_2}{\Delta}\Big]\\
 &=&\frac{1}{S'(0)}\left[-f'(t_1)-a'(0)f(t_1)\right]+O(\Delta).
\eea
Since $t_1=A$, for small $\Delta$, we have $w_1^*\approx P_A$, where
\be
 P_{A}=\frac{1}{f(A)S'(0)}\Big(-f'(A)-a'(0)f(A)\Big)=\frac{1}{2\lambda f(A)}\Big(-f'(A)+\lambda f(A)\Big).
\label{eq:ar1-pa}
\ee
Similarly, for the right boundary weight, we obtain
\bea
 {w_N^*}{}f(t_N)&=&\frac{f_N+k_1f_{N-1}}{S(\Delta)}=\frac{\Delta}{S(\Delta)}\frac{f_N-af_{N-1}}{\Delta}\\
 &=&\frac{\Delta}{S(\Delta)}\Big[\frac{f_N-f_{N-1}}{\Delta}+\frac{f_{N-1}-af_{N-1}}{\Delta}\Big]\\
 &=&
 \frac{1}{S'(0)}\left[f'(t_N)-a'(0)f(t_{N-1})\right]+O(\Delta).
\eea
Since $t_N=B$, for small $\Delta$, we have $w^*_N\approx P_B$, where
\be
 P_{B}=\frac{1}{f(B)S'(0)}\big(f'(B)-a'(0)f(B)\big)=\frac{1}{2\lambda f(B)}\big(f'(B)+\lambda f(B)\big).
\label{eq:ar1-pb}
\ee

Summarizing, we have proved the following result.

\begin{proposition} \label{prop1}
Consider the one-parameter regression model \eqref{eq:model-onepar} with $\mathrm{AR}(1)$ errors of the form \eqref{AR1},
where $0 < a <1$ and $f(\cdot)$ is a twice continuously differentiable function such that $f(t) \neq 0$ for all $t \in [A,B]$.  For large $N$,
the optimal discrete SLSE (defined in Lemma \ref{lem:2-1}) is approximated by
the continuous SLSE
\be
 \hat\theta_{}\!=\!D^*\Big(\!P_Af(A)y(A)\!+\!P_Bf(B)y(B)
 \!+\!\int_A^B \!\!p(t)f(t)y(t)dt
 \Big)
\label{eq:modif-slse-ar1}
\ee
where
\bea
 {D^*}\!=\! \Big(P_Af^2(A)\!+\!P_Bf^2(B)
  \!+\!\int_A^B p(t)f^2(t)dt\Big)^{-1},
\eea
and
$p(t)$, $P_A$ and $P_B$ are defined in
\eqref{eq:ar1-pt}, \eqref{eq:ar1-pa} and \eqref{eq:ar1-pb}, respectively.
For this  approximation, we have
$$D^*=\lim_{N\to\infty} \mathrm{Var}(\hat \theta_{\SLSE,N}),
$$
 i.e.
$D^*$ is the limit of the variance \eqref{eq:D-SLSw} of the optimal discrete SLSE design as $N\to\infty$.

\end{proposition}

\medskip

Throughout the following discussion  we call  a triple
$(p,P_A,P_B)$ containing a (signed) density $p$ and two weights $P_A$ and $P_B$, an approximate design for the  \emph{continuous}
SLSE estimator
defined in \eqref{eq:modif-slse-ar1}.

\medskip

\begin{remark}  \label{com} {\rm  Observing the discussion in the second part of Remark  \ref{arparameter}
it is reasonable
to compare  Proposition \ref{prop1}  with Theorem 2.1 in \cite{DetPZ2016}.
%\textbf{Remark 1 (connection with AOS paper \citep{DetPZ2016}).}
Note that the expressions for the optimal signed density $p(\cdot)$ and optimal weights $P_A$ and $P_B$ at boundary points
are  particular cases of the general formulae
\bea
 p(t)=-\frac{1}{  f(t ) v( t )} \Big[ \frac{h^{\prime }(t)  }{q^{\prime } (t)}  \Big]^\prime\, ,
\eea
\bea
P_A
= \frac{  1 }{ f(A)v^2(A) q^\prime (A)}  \Big[  \frac{f(A)u^\prime(A)}{u(A)} - f^\prime(A)  \Big]\, ,\;\;
P_B=   \frac{   h^\prime(B)}{f(B) v(B) q^\prime (B)}
\eea
with $u(t)=e^{\lambda t}$ and $v(s)=e^{-\lambda s}$,
where $q(t)=u(t)/v(t)$ and $h(t)=f(t)/v(t)$.
Indeed, we easily see that $h(t)=f(t)e^{\lambda t}$, $h'(t)=f'(t)e^{\lambda t}+f(t)\lambda e^{\lambda t}$,
$q'(t)=2\lambda e^{2\lambda t}$,
$h'(t)/q'(t)=f'(t)e^{-\lambda t}+f(t)\lambda e^{-\lambda t}$
and, consequently,
\bea
 p(t)&=&-\frac{1}{  f(t ) e^{-\lambda t}}(f'(t)e^{-\lambda t}+f(t)\lambda e^{-\lambda t})'\\
 &=&-\frac{1}{  f(t ) e^{-\lambda t}}(f''(t)e^{-\lambda t}-\lambda f''(t)e^{-\lambda t}+f'(t)\lambda e^{-\lambda t}-f(t)\lambda^2 e^{-\lambda t})\\
 &=&-\frac{1}{2\lambda f(t)}\Big(f''(t)-\lambda^2f(t)\Big).
\eea
as desired.
Similarly, we have
\bea
 P_A&=&\frac{  1 }{ f(A)e^{-2\lambda A} 2\lambda e^{2\lambda A}}  \Big[  \frac{f(A)\lambda e^{\lambda A}}{e^{\lambda A}} - f^\prime(A)  \Big]
 =\frac{1}{2\lambda f(A)}\Big(-f'(A)+\lambda f(A)\Big)
 \\
 P_B&=&   \frac{  f'(B)e^{\lambda B}+f(B)\lambda e^{\lambda B}}{f(B) e^{-\lambda B} 2\lambda e^{2\lambda B}}
 =\frac{1}{2\lambda f(B)}\Big(f'(B)+\lambda f(B)\Big)\, .
\eea
}
\end{remark}

%\newpage
\section{Autoregressive errors of order two}
\label{sec4}
\def\theequation{4.\arabic{equation}}
\setcounter{equation}{0}

In this section  we assume that the observations in   model \eqref{one} are  taken  at $N$  equidistant points  of the form
\eqref{points} and  that the errors  $\epsilon_1, \ldots , \epsilon_N$ satisfy
the discrete AR(2) equation
\be
\label{AR2}
 \epsilon_j-a_1\epsilon_{j-1}-a_2\epsilon_{j-2}=z_j,
\ee
%for all $j$ with artificial $\epsilon_0,\epsilon_{-1},\ldots$,
where $z_j$ are  Gaussian  independent identically distributed random variables
with mean $0$ and variance $\sigma^2_z=\sigma^2(1+a_2)((1-a^2)-a_1^2)/(1-a_2)$. Here we make a usual assumption  that \eqref{AR2} defines the AR(2) process for
$j \in \{ \ldots,-2,-1,0,1,2, \ldots\}$ but we only take the values such that $j \in \{ 1,2, \ldots, N\}$.
Let $r_k=\mathbb{E}[\epsilon_j\epsilon_{j+k}] $  be the autocovariance function of the AR$(2)$  process $\{ \epsilon_1,\ldots,\epsilon_N \}$
and assume without loss of generality that $\sigma^2=1$. It is well known that
the inverse of the covariance matrix $\Sigma =(\mathbb{E}[\epsilon_j\epsilon_{j}] )_{j,k} $  of the discrete AR$(2)$
 process is a  five-diagonal matrix, i.e.
\be
\label{inverseAR2}
 \Sigma^{-1}=\frac{1}{S}
 \begin{pmatrix}
 k_{11}&k_{12}&k_2&0  &0&0&\ldots\\
 k_{21}&k_{22}&k_1&k_2&0&0&\ldots\\
 k_2 & k_1    &k_0&k_1&k_2&0&\ldots\\
 0&k_2 & k_1    &k_0&k_1&k_2\\
 \vdots&\ddots&\ddots&\ddots&\ddots&\ddots&\ddots\\
 &&0&k_2& k_1    &k_0&k_1&k_2\\
 &&0&0&k_2&k_1&k_{22}&k_{12}\\
 &&0&0&0& k_2&k_{21}&k_{11}\\
 \end{pmatrix}~,
\ee
where the non-vanishing elements are given by $k_0=1+a_1^2+a_2^2$, $k_1=-a_1+a_1a_2$, $k_2=-a_2$,
$k_{11}=1$, $k_{12}=k_{21}=-a_1$, $k_{22}=1+a_1^2$ and
$S=(1+a_1-a_2)(1-a_1-a_2)(1+a_2)/(1-a_2)$.
Using Lemma \ref{lem:2-1} and the explicit form \eqref{inverseAR2} for  $\Sigma^{-1}$
we immediately obtain the following result.

\begin{corollary}
\label{cor:2.1}
Consider the linear regression  model \eqref{one} with  observations at $N$  equidistant points \eqref{points} and errors that follow
the discrete AR$(2)$ model  \eqref{AR2}.
  If $f_i =f(t_i) \neq 0 $ for $i=1,\dots,N$,  then  the  optimal weights in \eqref{eq:weights}   can be represented explicitly as follows:
\bea
%\label{w1}
 ~~~w_1^* &\!\!=\!\!&\frac{1}{Sf_1}  \left( k_{11} {f_1} + k_{12} {f_2} +k_2 f_3 \right)  \, , \\
%\label{w2}
 ~~~w_2^* &\!\!=\!\!&\frac{1}{Sf_2}  \left( k_{21} {f_1} + k_{22} {f_2} +k_1 f_3+k_2 f_4 \right) \, , \\
% \label{wN}
 ~~~w_N^* &\!\!=\!\!&\frac{1}{Sf_N}  \left( k_{11} {f_N} + k_{21} {f_{N-1}} +k_2 f_{N-2} \right)  \, ,\!\!\!\!\!\! \\
% \label{wN-1}
 ~~~w_{N-1}^* &\!\!=\!\!&\frac{1}{Sf_{N-1}}  \left( k_{12} {f_N} + k_{22} {f_{N-1}} +k_1 f_{N-2}+k_4 f_{N-3} \right)  \, ,\!\!\!\!\!\! \\
% \label{w30}
 ~~~w_i^*&\!\!=\!\!&\frac{1}{Sf_i}  \left(k_2 f_{i-2}+ k_1 f_{i-1}  +k_0 f_{i} + k_1 f_{i+1} + k_2 f_{i+2}\right) \nonumber
  \eea
  for $i=3, \ldots, N-2$.
\end{corollary}

For the approximation of $w^*_i$, we have
to study the behavior of  the coefficients which depend  on the autocovariance function $r_k$ of the
AR$(2)$ process \eqref{AR2}. There are different types of autocovariance functions  which will be introduced and
discussed in the remaining
part of this section. \\
Formally, a continuous AR$(2)$ process is a solution of the linear stochastic differential equation of the form
$$d \varepsilon'(t)=\tilde a_1\varepsilon'(t)+\tilde a_2\varepsilon(t)+\sigma^2_0dW(t),$$
 where $W(t)$ is a standard Wiener process,
[see \cite{BDY2007}].
Note that the process $\varepsilon(t)$ has the continuous derivative $\varepsilon'(t)$ and
the continuous process with drift is again denoted by $y(t) = \theta f(t) + \varepsilon (t) $,  ${t  \in [A,B]}$.
We also note that $y(t)$ is differentiable on the interval $ [A,B]$.

There are in fact three different  forms of the  autocovarince functions  (note that we assume throughout $\sigma^2=1$)
of  continuous AR(2) processes [see e.g. formulas (14)--(16) in \cite{he1989embedding}], which are given by
\be
 \rho^{(1)}(t)=\frac{\lambda_2}{\lambda_2-\lambda_1} e^{-\lambda_1|t|}-\frac{\lambda_1}{\lambda_2-\lambda_1} e^{-\lambda_2|t|}\, ,
\label{eq:KforAR2C}
\ee
 where $\lambda_1\neq\lambda_2$, $\lambda_1>0$, $\lambda_2>0$, by
\bea
\rho^{(2)}(t)=e^{-\lambda |t|}\Big\{\cos(q |t|)+ \frac{\lambda}{q} \sin(q |t|)\Big\}\,,
%\label{eq:K-ecosC}
\eea
where $\lambda>0$, $q>0$, and by
\bea
\rho^{(3)}(t)=e^{-\lambda |t|}(1+ \lambda |t|)\, ,
%\label{eq:K-elinC}
\eea
where $\lambda>0$. From   formulas (11)--(13) in \cite{he1989embedding} we obtain that
the   corresponding three forms of the autocovariances
 of the discrete  AR(2) process  of the form \eqref{AR2} are given by
\be
 r_k^{(1)}=  \mathbb{E}[\epsilon_j\epsilon_{j+k}]  = C p_1^k+(1-C) p_2^k, ~~~C=\frac{(1-p_2^2)p_1}{(1-p_2^2)p_1-(1-p_1^2)p_2}\, ,
\label{eq:KforAR2}
\ee
where $j\geq 0$, $p_1\neq p_2$, $0<|p_1|,|p_2|<1$; by
\be
r^{(2)}_k=p^k\big(\cos(bk)+C\sin(bk)\big),~~~ C=\cot(b)\frac{1-p^2}{1+p^2} \, ,
\label{eq:K-ecos}
\ee
where $0<p<1$, $0<b<2\pi$ and $b \neq \pi$, and finally by
\be
r^{(3)}_k=p^k\,(1+k C),~~~ C=\frac{1-p^2}{1+p^2}\,,
\label{eq:K-elin}
\ee
where $0<|p|<1$. In the following subsections we determine approximations for the optimal weights $w_i^*$  in  Lemma \ref{lem:2-1}
for the different types of autocovariance functions. All results will be summarized in Theorem \ref{th:approx-ar2} below.

\subsection{Autocovariances of the form \eqref{eq:KforAR2}}

From Corollary \ref{cor:2.1} we obtain that
\bea
 S{w_i^*}{}f_i&=&-a_2f_{i-2}+(a_1a_2-a_1)f_{i-1}+(1+a_1^2+a_2^2)f_i+(a_1a_2-a_1)f_{i+1}-a_2f_{i+2}\\
 &=&a_2(2f_i-f_{i-2}-f_{i+2})-(a_1a_2-a_1)(2f_i-f_{i-1}-f_{i+1})\\
 &&+(1+a_1^2+a_2^2-2a_2+2a_1a_2-2a_1)f_i\\
 &=&a_2(2f_i-f_{i-2}-f_{i+2})-(a_1a_2-a_1)(2f_i-f_{i-1}-f_{i+1})\\
 &&+(a_1+a_2-1)^2f_i
\eea
for $i=3,4,\ldots,N-2.$
Now consider the case when
the autocovariance structure of the errors has the form \eqref{eq:KforAR2} for fixed $N$.
Suppose that the parameters of the autocovariance function \eqref{eq:KforAR2} satisfy
$p_1\neq p_2$, $0<p_1,p_2<1$.
We do not discuss the case with negative $p_1$ or negative $p_2$ because discrete AR$(2)$
 processes with such parameters
do not have   continuous real-valued analogues.
From the Yule-Walker equations we obtain that the coefficients $a_1$ and $a_2$ in \eqref{AR2}
are given by
\be
  a_1=r_1\frac{1-r_2}{1-r^2_1},~~~ a_2=\frac{r_2-r^2_1}{1-r^2_1}\, ,
  \label{eq:a1,a2}
\ee
where $r_1=r_1^{(1)}$ and $r_2=r_2^{(1)}$ are defined  by \eqref{eq:KforAR2}.
With the notation
$\lambda_1= -\log(p_1) / \Delta$ and $\lambda_2= -\log(p_2) / \Delta$ with $\Delta=(B-A)/N$
we obtain
\be
\label{T0}
 p_1=e^{-\lambda_1\Delta}, ~~~p_2=e^{-\lambda_2\Delta}.
\ee
We will assume that $\lambda_1$ and $\lambda_2$ are fixed but $\Delta$ is small and consider the properties of different quantities as  $\Delta \to 0$.
By a straightforward   Taylor expansion  we obtain the approximations
\be
%\label{T1}
 a_1=a_1(\Delta)&=&2-(\lambda_1+\lambda_2)\Delta+(\lambda_1^2+\lambda_2^2)\Delta^2/2+O(\Delta^3),\nonumber \\
%\label{T2}
a_2=a_2(\Delta)&=&-1+(\lambda_1+\lambda_2)\Delta-(\lambda_1+\lambda_2)^2\Delta^2/2+O(\Delta^3), \nonumber\\
% \label{T3}
 S=S(\Delta)&=&2\lambda_1\lambda_2(\lambda_1+\lambda_2)\Delta^3+O(\Delta^4), \nonumber\\
\label{T4} C=C(\Delta)&=&\frac{\lambda_2}{\lambda_2-\lambda_1}+\frac{1}{6}\lambda_1\lambda_2\frac{\lambda_1+\lambda_2}{\lambda_1-\lambda_2}\Delta^2+O(\Delta^4).
\ee
Consequently  (observing  \eqref{T0} and \eqref{T4}),  for large $N$ the continuous AR(2) process with autocovariances \eqref{eq:KforAR2C} can be considered as an
approximation to the discrete AR(2) process with autocovariances
\eqref{eq:KforAR2}.

Since $S=O(\Delta^3)$, $a_1=2+O(\Delta)$ and $a_2=-1+O(\Delta)$, % and dividing the previous equality on $\Delta^4$,
it follows
\bea
 S\frac{w_i^*}{\Delta^4}f_i
 &=&-4a_2\frac{1}{\Delta^2}f''(t_i)+(a_1a_2-a_1)\frac{1}{\Delta^2}f''(t_i)+\frac{1}{\Delta^4}(a_1+a_2-1)^2f_i+O(\Delta)\\
 &=&\frac{1}{\Delta^2}(a_1a_2-a_1-4a_2)f''(t_i)+\frac{1}{\Delta^4}(a_1+a_2-1)^2f_i+O(\Delta)\\
 &=&-(\lambda_1^2+\lambda_2^2)f''(t_i)+\lambda_1^2\lambda_2^2f_i+O(\Delta).
\eea
Thus, the optimal weights $w_i^*$, $i=3,\ldots,N-2,$ are approximated by the signed density
\be
 p(t)=-\frac{1}{s_3 f(t)}\big( (\lambda_1^2+\lambda_2^2)f''(t)-\lambda_1^2\lambda_2^2 f(t)\big),
 \label{eq:ar2-pt}
\ee
where $s_3=2\lambda_1\lambda_2(\lambda_1+\lambda_2)$.
For  the boundary points we obtain
\bea
 S{w_1^*}{}f_1&=&f_1-a_1f_2-a_2f_3\\
 &=&(-2f_2+f_3+f_1)+(\lambda_1+\lambda_2)(f_2-f_3)\Delta\\
 &&+((-1/2 f_2+1/2 f_3) \lambda_1^2+f_3 \lambda_1 \lambda_2+(-1/2 f_2+1/2 f_3) \lambda_2^2)\Delta^2\\
 &&+((1/6 f_2-1/6 f_3) \lambda_1^3-1/2 f_3 \lambda_1^2 \lambda_2-1/2 f_3 \lambda_1 \lambda_2^2+(1/6 f_2-1/6 f_3) \lambda_2^3)\Delta^3
 +O(\Delta^4)\\
 &=&\big(f''(t_2)-(\lambda_1+\lambda_2)f'(t_2)+f_3 \lambda_1 \lambda_2\big) \Delta^2+O(\Delta^3)
\eea
and
\bea
 S{w_2^*}{}f_2&=&-a_1f_1+(1+a_1^2)f_2+(a_1a_2-a_1)f_3-a_2f_4\\
 &=&(-2 f_1+f_4+5 f_2-4 f_3)+(\lambda_1+\lambda_2) (f_1-4 f_2+4 f_3-f_4) \Delta\\
&&+((-1/2 f_1+1/2 f_4-3 f_3+3 f_2) \lambda_1^2+(2 f_2-4 f_3+f_4) \lambda_2 \lambda_1\\
&&~~~ +(-1/2 f_1+1/2 f_4-3 f_3+3 f_2) \lambda_2^2) \Delta^2\\
&&+((1/6 f_1-5/3 f_2+5/3 f_3-1/6 f_4) \lambda_1^3+(-f_2+3 f_3-1/2 f_4) \lambda_2 \lambda_1^2\\
&&~~~ +(-f_2+3 f_3-1/2 f_4) \lambda_2^2 \lambda_1+(1/6 f_1-5/3 f_2+5/3 f_3-1/6 f_4) \lambda_2^3) \Delta^3 +O(\Delta^4)\\
&=&\big(f''(t_3)-2f''(t_2)+ (\lambda_1+\lambda_2)(3f'(t_2)-f'(t_1)-f'(t_3))-f_3 \lambda_1 \lambda_2 \big)\Delta^2+O(\Delta^3)\\
&=&\big(-f''(t_2)+ (\lambda_1+\lambda_2)f'(t_2)-f_3 \lambda_1 \lambda_2 \big)\Delta^2+O(\Delta^3)
\eea
Thus, we can see that
$$
{w_1^*}=- {w_2^*}+O(1)= Q_A\frac{1}{\Delta}+O(1),
$$
where
\be
 Q_A=\frac{1}{s_3f(A)}\big(f''(A)- (\lambda_1+\lambda_2)f'(A)+ \lambda_1 \lambda_2f(A)\big).
 \label{eq:qa-expexp}
\ee
This means that the coefficients $w_1^*$ and $w_2^*$ at $t_1$ and $t_2$ are large in absolute value and have different signs.
Similarly, we have
$$
w_N^* =- w_{N-1}^*+O(1) = Q_B\frac{1}{\Delta}+O(1)
$$
where
\be
 Q_B=\frac{1}{s_3f(B)}\big(f''(B)+ (\lambda_1+\lambda_2)f'(B)+ \lambda_1 \lambda_2f(B)\big).
 \label{eq:qb-expexp}
\ee

%We can approximate these BLUE coefficients by other robust coefficients.
%We can do this because there is no jump in the process $y(t)$.
To do a finer approximation, we have to investigate the quantity
\bea
g:=S{w_1^*}{}f_1+S{w_2^*}{}f_2,
\eea
which is  of order $O(1)$.
Indeed, we have
\bea
g&=&(3 f_2-3 f_3-f_1+f_4)+(\lambda_2+\lambda_1) (f_1-3 f_2+3 f_3-f_4) \Delta\\
&&+((- f_1+ f_4-5 f_3+5 f_2)/2 (\lambda_1^2+\lambda_2^2)+(2 f_2-3 f_3+f_4) \lambda_2 \lambda_1) \Delta^2\\
&&+(( f_1-9 f_2+9 f_3- f_4)/6 (\lambda_1^3+\lambda_2^3)\\
&&+(-2 f_2+5 f_3- f_4)/2 (\lambda_1^2 \lambda_2+\lambda_1 \lambda_2^2) ) \Delta^3
+O(\Delta^4)\\
&=&f'''(t_1)\Delta^3+O(\Delta^4)
+(-f'(t_1) (\lambda_1^2+\lambda_2^2)-f'(t_1) \lambda_2 \lambda_1) \Delta^3 \\
&&+f(t_1)(\lambda_1^2 \lambda_2+\lambda_1 \lambda_2^2)  \Delta^3
+O(\Delta^4)\\
&=& \big(f'''(t_1) -(\lambda_1^2+\lambda_1\lambda_2+\lambda_2^2) f'(t_1)  + \lambda_1 \lambda_2(\lambda_1+ \lambda_2)f(t_1)\big)\Delta^3+O(\Delta^4)
\eea
and, consequently,
\bea
 {w_1^*}{}f_1+{w_2^*}{}f_2=
 \frac{1}{s_3}\big(f'''(t_1) -(\lambda_1^2+\lambda_1\lambda_2+\lambda_2^2) f'(t_1)  + \lambda_1 \lambda_2(\lambda_1+ \lambda_2)f(t_1)\big)+O(\Delta),
\eea
where $s_3=2\lambda_1\lambda_2(\lambda_1+\lambda_2)$.
Therefore, if $\Delta \to 0$, it follows that $w_1^*+w_2^*\approx P_A$, where
\be
 P_{A}=\frac{1}{s_3f(A)}\big(f'''(A) -(\lambda_1^2+\lambda_1\lambda_2+\lambda_2^2) f'(A)  + \lambda_1 \lambda_2(\lambda_1+ \lambda_2)f(A)\big).
 \label{eq:ar2-pa}
\ee
Similarly, we obtain
  $w_N^*+w_{N-1}^*\approx P_B$ if $\Delta \to 0$, where
\be
 P_{B}=\frac{1}{s_3f(B)}\big(-f'''(B) +(\lambda_1^2+\lambda_1\lambda_2+\lambda_2^2) f'(B)  + \lambda_1 \lambda_2(\lambda_1+ \lambda_2)f(B)\big).
 \label{eq:ar2-pb}
\ee
Summarizing, we have proved the following result.

\medskip

\begin{proposition} \label{prop2}
Consider the one-parameter model \eqref{one} such that the errors follow the  $\mathrm{AR}(2)$ model
with  autocovariance function  \eqref{eq:KforAR2}.
Assume that $f(\cdot)$ is a three times  continuously differentiable and  $f(t) \neq 0$ for all $t \in [A,B]$.
Then for large $N$,
the optimal discrete SLSE (defined in Lemma \ref{lem:2-1}) can be  approximated by
the  continuous SLSE
{\small
\be
 \hat\theta_{}\!=\!D^*\Big(Q_Bf(B)y'(B)\!-\!Q_Af(A)y'(A)\!+\!P_Af(A)y(A)\!+\!P_Bf(B)y(B)
 \!+\!\int_A^B \!\!\!\!p(t)f(t)y(t)dt
 \Big) ~~~
\label{eq:modif-slse}
\ee}
where
\bea
 D^*\!=\! \Big(Q_Bf(B)f'(B)\!-\!Q_Af(A)f'(A)\!+\!P_Af^2(A)\!+\!P_Bf^2(B)
  \!+\!\int_A^B p(t)f^2(t)dt\Big)^{-1}
\eea
and
$p(t)$, $Q_A$, $Q_B$, $P_A$ and $P_B$  defined in
\eqref{eq:ar2-pt}, \eqref{eq:qa-expexp}, \eqref{eq:qb-expexp}, \eqref{eq:ar2-pa} and \eqref{eq:ar2-pb} respectively.
For this approximation, we have  $D^*=\lim_{N\to\infty} \mathrm{Var}(\hat \theta_{\SLSE,N})$, i.e.
$D^*$ is the limit of the variance \eqref{eq:D-SLSw} of the optimal discrete SLSE design as $N\to\infty$.
\end{proposition}

\medskip

In the following  discussion we call  a tuple
($p,Q_A, Q_B, P_A, P_B$)
which contains  a (signed) density $p(\cdot)$ and four  weights $Q_A, Q_B, P_A, P_B$, an approximate design for the  continuous
SLSE estimator
defined in
 \eqref{eq:modif-slse}.

\subsection{Autocovariances of the form \eqref{eq:K-ecos}}

Consider the autocovariance function of  the form \eqref{eq:K-ecos}, then the coefficients $a_1$ and $a_2$
are given by \eqref{eq:a1,a2} where
$r_1=r_1^{(2)} $ and $r_2=r_2^{(2)} $ are defined by \eqref{eq:K-ecos}.  With the notations
$\lambda= -\log p / \Delta$ and $q=b/ \Delta$ (or equivalently
 $p=e^{-\lambda\Delta}$ and $b=q\Delta$) we obtain by  a Taylor expansion
\bea
 a_1&=&2-2\lambda\Delta+(\lambda^2-q^2)\Delta^2 +O(\Delta^3),\\
 a_2&=&-1+2\lambda\Delta-2\lambda^2\Delta^2+O(\Delta^3),\\
 S&=&4\lambda(\lambda^2+q^2) \Delta^3+O(\Delta^4)
\eea
and
\bea
 C=\frac{\lambda}{q}- \frac{\lambda (\lambda^2+q^2)}{3q} \Delta^2+O(\Delta^4)
\eea
as $\Delta \to 0$.
Similarly, we have
\bea
 S\frac{w_i^*}{\Delta^4}f_i
 &=&\frac{1}{\Delta^2}(a_1a_2-a_1-4a_2)f''(t_i)+\frac{1}{\Delta^4}(a_1+a_2-1)^2f_i+O(\Delta)\\
 &=&-2(\lambda^2-q^2)f''(t_i)+(\lambda^2+q^2)^2 f_i+O(\Delta).
\eea
Thus, the optimal weights $w_i^*$, $i=3,\ldots,N-2,$ are approximated by the signed density
\be
 p(t)=-\frac{1}{s_3 f(t)}\big( 2(\lambda^2-q^2)f''(t)-(\lambda^2+q^2)^2 f(t)\big),
 \label{eq:ar2-pt-cos}
\ee
where $s_3=4\lambda(\lambda^2+q^2)$.
Similarly, we obtain that
\bea
{w_1^*}&=& - {w_2^*}+O(1)= Q_A\frac{1}{\Delta}+O(1), \\
{w_N^*}&=& - {w_{N-1}^*}+O(1)= Q_B\frac{1}{\Delta}+O(1),
\eea
where
\be
 Q_A &=&\frac{1}{s_3f(A)}\big(f''(A)- 2\lambda f'(A)+ (\lambda^2+q^2) f(A)\big),
 \label{eq:ar2-qa-cos}
\\
 Q_B&=&\frac{1}{s_3f(B)}\big(f''(B)+ 2\lambda f'(B)+ (\lambda^2+q^2) f(B)\big).
 \label{eq:ar2-qb-cos}
\ee

Calculating $g:=S{w_1^*}{}f_1+S{w_2^*}{}f_2$ we have
\bea
 g&=&(3f_2-f_1-3f_3+f_4)+2\lambda(f_1-3f_2+3f_3-f_4)\Delta\\
 &&+((-f_1+7 f_2-8 f_3+2f_4)\lambda^2+q^2(f_1-3 f_2+2f_3))\Delta^2\\
 &&+((-f_1+7 f_2-4 f_3)\lambda q^2+( f_1-15 f_2+24 f_3-4f_4)/3 \lambda^3)\Delta^3+O(\Delta^4)\\
 &=&f'''(t_1)\Delta^3- (3\lambda^2-q^2)f'(t_1)\Delta^3+2\lambda(\lambda^2+q^2)f(t_1)\Delta^3+O(\Delta^4).
\eea

Therefore,   it follows that  $w_1^*+w_2^*\approx P_A$ if $\Delta \to 0$, where
\be
 P_{A}=\frac{1}{s_3f(A)}\big(f'''(A) -(3\lambda^2-q^2) f'(A)  + 2\lambda(\lambda^2+q^2) f(A)\big),
 \label{eq:ar2-pa-cos}
\ee
and  $s_3=4\lambda(\lambda^2+q^2)$.
Similarly, we obtain
  the approximation $w_N^*+w_{N-1}^*\approx P_B$ if $\Delta \to 0$, where
\be
 P_{B}=\frac{1}{s_3f(B)}\big(-f'''(B) +(3\lambda^2-q^2) f'(B)  + 2\lambda(\lambda^2+q^2) f(B)\big).
 \label{eq:ar2-pb-cos}
\ee
Summarizing, we have proved the following result.

\medskip

\begin{proposition}  \label{prop3}
Consider the one-parameter model \eqref{one} such that the errors follow the  $\mathrm{AR}(2)$ model
with  autocovariance function  \eqref{eq:K-ecos}.
Assume that $f(\cdot)$ is a three times  continuously differentiable and  $f(t) \neq 0$ for all $t \in [A,B]$.
Then for large $N$,
the optimal discrete SLSE (defined in Lemma \ref{lem:2-1}) can be approximated by
the continuous SLSE \eqref{eq:modif-slse}, where the tuple ($p, Q_A, Q_B, P_A, P_B$)
is  defined by
\eqref{eq:ar2-pt-cos}, \eqref{eq:ar2-qa-cos}, \eqref{eq:ar2-qb-cos}, \eqref{eq:ar2-pa-cos} and \eqref{eq:ar2-pb-cos} respectively.

\end{proposition}

\medskip

\subsection{Autocovariances of the form \eqref{eq:K-elin}}

For the autocovariance function \eqref{eq:K-elin}
the coefficients $a_1$ and $a_2$
in the AR$(2)$ process are given by \eqref{eq:a1,a2} where
$r_1=r_1^{(3)} $ and $r_2=r_2^{(3)}$ are defined by \eqref{eq:K-elin}.
With the notation  $\lambda =- \log p / \Delta$ (or equivalently  $p=e^{-\lambda\Delta}$)
we obtain the Taylor expansions
\bea
 a_1&=&2-2\lambda\Delta+\lambda^2\Delta^2+O(\Delta^3),\\
 a_2&=&-1+2\lambda\Delta-2\lambda^2\Delta^2+O(\Delta^3),\\
 S&=&4\lambda^3\Delta^3+O(\Delta^4),\\
 C&=&\lambda\Delta-\frac{\lambda^3}{3}\Delta^3+O(\Delta^5)
\eea
as $\Delta \to 0$. Similar calculations as given in the previous paragraphs give
\bea
 S\frac{w_i^*}{\Delta^4}f_i
 &=&\frac{1}{\Delta^2}(a_1a_2-a_1-4a_2)f''(t_i)+\frac{1}{\Delta^4}(a_1+a_2-1)^2f_i+O(\Delta)\\
 &=&-2\lambda^2f''(t_i)+\lambda^4 f_i+O(\Delta).
\eea
Thus, the optimal weights $w_i^*$, $i=3,\ldots,N-2,$ are approximated by the signed density
\be
 p(t)=-\frac{1}{s_3 f(t)}\big( 2\lambda^2f''(t)-\lambda^4 f(t)\big),
 \label{eq:ar2-pt-lin}
\ee
where $s_3=4\lambda^3$. For the remaining weighs $w_1^*$ and $w_2^*$  we obtain
\bea
{w_1^*}&=& - {w_2^*}+O(1)= Q_A\frac{1}{\Delta}+O(1), \\
{w_N^*} &=& - {w_{N-1}^*}+O(1)= Q_B\frac{1}{\Delta}+O(1),
\eea
with
\be
 Q_A &=&\frac{1}{s_3f(A)}\big(f''(A)- 2\lambda f'(A)+ \lambda^2 f(A)\big),
 \label{eq:ar2-qa-lin} \\
 Q_B&=&\frac{1}{s_3f(B)}\big(f''(B)+ 2\lambda f'(B)+ \lambda^2 f(B)\big).
 \label{eq:ar2-qb-lin}
\ee
Calculating $g:=S{w_1^*}f_1+S{w_2^*}f_2$ we have
\bea
 g&=&(3f_2-3f_3-f_1+f_4)+2\lambda(f_1-3f_2+3f_3-f_4)\Delta\\
 &&-\lambda^2(f_1-7f_2+8f_3-2f_4)\Delta^2\\
 &&+1/3\lambda^3(f_1-15f_2+24f_3-4f_4)\Delta^3 +O(\Delta^4)\\
 &=&f'''(t_1)\Delta^3-3\lambda^2f'(t_1)\Delta^3+2\lambda^3f(t_1)\Delta^3 +O(\Delta^4).
\eea
Therefore, if $\Delta \to 0$, it follows that  $w_1^*+w_2^*\approx P_A$, where
\be
 P_{A}=\frac{1}{s_3f(A)}\big(f'''(A) -3\lambda^2 f'(A)  + 2\lambda^3 f(A)\big),
 \label{eq:ar2-pa-lin}
\ee
and   $s_3=4\lambda^3$.
Similarly, we obtain
  the approximation $w_N^*+w_{N-1}^*\approx P_B$ if $\Delta \to 0$, where
\be
 P_{B}=\frac{1}{s_3f(B)}\big(-f'''(B) +3\lambda^2 f'(B)  + 2\lambda^3 f(B)\big).
 \label{eq:ar2-pb-lin}
\ee
Summarizing, we have proved the following result.

\begin{proposition}  \label{prop4}
Consider the one-parameter model \eqref{one} such that the errors follow the  $\mathrm{AR}(2)$ model
with  autocovariance function  \eqref{eq:K-elin}.
Then for large $N$,
the optimal discrete SLSE (defined in Lemma \ref{lem:2-1}) can be approximated by
the  continuous SLSE \eqref{eq:modif-slse},
where the tuple $(p, Q_A, Q_B, P_A, P_B$) is  defined in
\eqref{eq:ar2-pt-lin}, \eqref{eq:ar2-qa-lin}, \eqref{eq:ar2-qb-lin}, \eqref{eq:ar2-pa-lin} and \eqref{eq:ar2-pb-lin} respectively.
\end{proposition}

\subsection{General statement}

Propositions \ref{prop2} - \ref{prop4}   can be combined in  the following statement.

\begin{theorem}
\label{th:approx-ar2}
Consider the one-parameter model \eqref{one} such that the errors follow the  $\mathrm{AR}(2)$ model.
Assume that $f(\cdot)$ is a three times  continuously differentiable and  $f(t) \neq 0$ for all $t \in [A,B]$.
Define the following constants depending on the form of the autocovariance function $r_k$.
If $r_k$ is of the form \eqref{eq:KforAR2}, set
\bea
 \lambda_1 &=&  -\frac{\ln(p_1)}{\Delta},~\lambda_2= -\frac{\ln(p_2)}{\Delta}, \\
 \tau_0 &=& \lambda_1^2\lambda_2^2,~\tau_2=\lambda_1^2+\lambda_2^2,~\beta_1=\lambda_1+\lambda_2,~\beta_0=\lambda_1\lambda_2,\\
 \gamma_1&=& \lambda_1^2+\lambda_1\lambda_2+\lambda_2^2~,\gamma_0=\lambda_1 \lambda_2(\lambda_1+ \lambda_2),~
 s_3=2\lambda_1\lambda_2(\lambda_1+\lambda_2).
\eea
If $r_k$ is of  the form \eqref{eq:K-ecos}, set
\bea
 \lambda &=& -\frac{\ln(p)}{\Delta},~q= -\frac{b}{\Delta}, \\
 \tau_0&=& (\lambda^2+q^2)^2,~\tau_2=2(\lambda^2-q^2),~\beta_1=2\lambda,~\beta_0=\lambda^2+q^2,\\
 \gamma_1&=& (3\lambda^2-q^2)~,\gamma_0=2\lambda(\lambda^2+q^2),~
 s_3=4\lambda(\lambda^2+q^2).
\eea
If $r_k$ is of the form \eqref{eq:K-elin}, set
\bea
 \lambda &=& -\frac{\ln(p)}{\Delta},~
 \tau_0=\lambda^4,~\tau_2=2\lambda^2,~\beta_1=2\lambda,~\beta_0=\lambda^2,\\
 \gamma_1&=&3\lambda^2~,\gamma_0=2\lambda^3,~
 s_3=4\lambda^3.
\eea
For large $N$,
the optimal discrete SLSE (defined in Lemma \ref{lem:2-1}) can be approximated by
the continuous SLSE
{\small
\bea
 \hat\theta_{}\!=\!D^*\big(Q_Bf(B)y'(B)\!-\!Q_Af(A)y'(A)\!+\!P_Af(A)y(A)\!+\!P_Bf(B)y(B)
 \!+\!\int_A^B \!\!\!\!p(t)f(t)y(t)dt
 \big)
%\label{eq:modif-slse-gen}
\eea}
where
\bea
 D^*\!=\! \Big(Q_Bf(B)f'(B)\!-\!Q_Af(A)f'(A)\!+\!P_Af^2(A)\!+\!P_Bf^2(B)
  \!+\!\int_A^B p(t)f^2(t)dt\Big)^{-1}.
\eea
For this approximation, we have  $D^*=\lim_{N\to\infty} \mathrm{Var}(\hat \theta_{SLSE,N})$, i.e.
$D^*$ is the limit of the variance \eqref{eq:D-SLSw} of the optimal discrete SLSE design as $N\to\infty$.
Here  the quantities
$p(t)$, $Q_A$, $Q_B$, $P_A$ and $P_B$ in the continuous SLSE are  defined by
\be
 p(t) &=& -\frac{1}{s_3 f(t)}\big( \tau_2 f''(t)-\tau_0 f(t)\big), \\
 \label{eq:ar2-ptgen}
 P_{A} &=& \frac{1}{s_3f(A)}\big(f'''(A) -\gamma_1 f'(A)  + \gamma_0 f(A)\big), \nonumber \\
 P_{B}&=&\frac{1}{s_3f(B)}\big(-f'''(B) +\gamma_1 f'(B)  + \gamma_0 f(B)\big), \nonumber \\
 Q_A&=& \frac{1}{s_3f(A)}\big(f''(A)- \beta_1f'(A)+ \beta_0 f(A)\big), \nonumber \\
 Q_B &=&\frac{1}{s_3f(B)}\big(f''(B)+ \beta_1f'(B)+ \beta_0 f(B)\big).
\ee
\end{theorem}

\section{Examples} \label{sec5}

\subsection{Approximations of the discrete SLSE }

Consider the one-parameter model with $f(t)=t^{\alpha}$ and AR(1) errors $(0 < a < 1)$.
The design space is given by an interval $[A,B]$ such that $f(t)\neq0$ for all $t\in[A,B]$.
Then the optimal discrete design for the SLSE is approximated by
a  design of the form
\eqref{cont_design},
where the density $p(t)$,  and the weights $P_A$ and $P_B$ are defined by
\bea
 p(t)&=&-\frac{1}{2\lambda }\big(  \alpha(\alpha-1)t^{-2}-\lambda^2 \big),\\
 P_{A}&=&\frac{1}{2\lambda}\big( - \alpha A^{-1}  + \lambda \big),\\
 P_{B}&=&\frac{1}{2\lambda}\big(  \alpha B^{-1}  + \lambda \big).
\eea
In Table~\ref{tab:optdes-ar1} we display values of $p(t)$, $P(A)$ and $P_B$ for several
exponents $\alpha$ and also for  the regression function $f(t)=e^t$
For example, if $f(t)=e^t$ we  observe that  $P_A$ is positive for $\lambda>1$ and negative for $0<\lambda<1$,
$P_B$ is positive for $\lambda>0$,
$p(t)$ is positive for $\lambda>1$ and negative for $\lambda\in(0,1)$.
For large $\lambda$,  the contribution of observations
at the interval $(A,B)$ to the continuous  SLSE is significant.
For the location model $f(t)=1$, we can see that $P_B=P_B=1/2$ and $p(t)=\lambda/2$.
This implies that for small $\lambda$ the contribution of observations
at boundary points to the continuous SLSE is large
and the contribution of observations
at the interval $(A,B)$ to the continuous  SLSE is small.
For large $\lambda$,  the contribution of observations
at the interval $(A,B)$ to the continuous  SLSE is essential.
\begin{table}[!hhh]
\caption{\it The function $p(t)$ and the weights $P_A$ and $P_B$ of the continuous   SLSE  for several functions $f(t)$ and an AR(1) error process.}
%with the autocorrelation function in the form \eqref{eq:KforAR2},
\medskip
\centering
\begin{tabular}{cccc}
\hline
$f(t)$&$P_A$&$P_B$&$p(t)$\\
\hline
1&$\ds\frac12$&$\ds\frac12$&$\ds\frac{\lambda}{2}$\rule{0mm}{8mm}\\
$t$&$\ds\frac12-\frac{1}{2A\lambda}$&$\ds\frac12+\frac{1}{2B\lambda}$&$\ds\frac{\lambda}{2}$\rule{0mm}{8mm}\\
$t^2$&$\ds\frac12-\frac{1}{A\lambda}$&$\ds\frac12+\frac{1}{B\lambda}$&$\ds\frac{\lambda}{2}-\frac{1}{\lambda t^2}$\rule{0mm}{8mm}\\
$t^3$&$\ds\frac12-\frac{3}{2A\lambda}$&$~~~\ds\frac12+\frac{3}{2B\lambda}~~~$&$\ds\frac{\lambda}{2}-\frac{3}{\lambda t^2}$\rule{0mm}{8mm}\\
$t^4$&$\ds\frac12-\frac{2}{A\lambda}$&$\ds\frac12+\frac{2}{B\lambda}$&$\ds\frac{\lambda}{2}-\frac{6}{\lambda t^2}$\rule{0mm}{8mm}\\
$~~e^t~~$&$~~\ds\frac12-\frac{1}{2\lambda}~~$&$~~\ds\frac12+\frac{1}{2\lambda}~~$&$~~\ds\frac{\lambda}{2}-\frac{1}{2\lambda }~~$\rule{0mm}{8mm}\\
\\
\hline
\end{tabular}
\label{tab:optdes-ar1}
\end{table}

Next we consider the same models  with  an AR$(2)$ error process.
For example, if $f(t)=t^{\alpha}$
the SLSE   is approximated by
the continuous SLSE of the form
\eqref{cont_design}, where
\bea
 p(t)&=&-\frac{1}{s_3 }\big( \tau_2 \alpha(\alpha-1)t^{-2}-\tau_0 \big),\\
 P_{A}&=&\frac{1}{s_3}\big(\alpha(\alpha-1)(\alpha-2)A^{-3} -\gamma_1 \alpha A^{-1}  + \gamma_0 \big),\\
 P_{B}&=&\frac{1}{s_3}\big(-\alpha(\alpha-1)(\alpha-2)B^{-3} +\gamma_1 \alpha B^{-1}  + \gamma_0 \big),\\
 Q_A&=&\frac{1}{s_3}\big(\alpha(\alpha-1)A^{-2} -\beta_1 \alpha A^{-1}  + \beta_0 \big),\\
 Q_B&=&\frac{1}{s_3}\big(\alpha(\alpha-1)B^{-2} +\beta_1 \alpha B^{-1}  + \beta_0 \big).\\
\eea
Note that signs of $p(t),$ $Q_A$, $Q_B$, $P_A$ and $P_B$ depend on the form of the autocovariance function and its parameters.
For  the form \eqref{eq:K-elin},
we provide values of $p(t)$, $Q_A$, $Q_B$, $P_A$ and $P_B$ for several functions $f(t)$
in Table~\ref{tab:optdes-ar2}. The other cases can be obtained similarly and are not displayed for the sake of brevity.

For example, if $f(t)=e^t$ we can see  that both $P_A$ and $Q_A$ are positive for all $\lambda\neq1$,
$P_B$ is positive for $\lambda>0.5$ and negative for $\lambda\in(0,0.5)$,
$p(t)$ is positive for $\lambda>\sqrt{2}$ and negative for $\lambda\in(0,\sqrt{2})$.
For large $\lambda$,  the contribution of observations
at the interval $(A,B)$ to the continuous  SLSE is notable.
For the location model $f(t)=1$, we can see that $P_A=P_B=1/2$, $Q_A=Q_B=1/(4\lambda)$ and $p(t)=\lambda/4$.
This implies that for small $\lambda$ the contribution of observations
at boundary points to the continuous  SLSE is very large
and the contribution of observations
at the interval $(A,B)$ to the continuous SLSE is small.
For large $\lambda$,  the contribution of observations
at the interval $(A,B)$ to the continuous SLSE is essential.

\begin{table}[!hhh]
\caption{\it The function $p(t)$ and the weights $P_A$, $P_B$, $Q_A$ and $Q_B$ in the continuous SLSE  for several functions $f(t)$ and an AR(2) error process with the autocovariance function \eqref{eq:K-elin}.}
%with the autocorrelation function in the form \eqref{eq:KforAR2},
\medskip

\centering
\begin{tabular}{cccccc}
\hline
$f(t)$&$P_A$&$P_B$&$p(t)$ & $Q_A$&$Q_B$\\
\hline
1&$ \frac12$&$ \frac12$&$ \frac{\lambda}{4}$\rule{0mm}{5mm}
&$\frac{1}{4\lambda}$&$\frac{1}{4\lambda}$\\
$t$&$ \frac12-\frac{3}{4A\lambda}$&$ \frac12+\frac{3}{4B\lambda}$&$ \frac{\lambda}{4}$\rule{0mm}{5mm}
&$\frac{1}{4\lambda}-\frac{1}{2A\lambda^2}$&$\frac{1}{4\lambda}+\frac{1}{2B\lambda^2}$\\
$t^2$&$ \frac12-\frac{3}{2A\lambda}$&$ \frac12+\frac{3}{2B\lambda}$&$ \frac{\lambda}{4}-\frac{1}{\lambda t^2}$\rule{0mm}{5mm}
&$\frac{1}{4\lambda}-\frac{1}{A\lambda^2}+\frac{1}{2A^2\lambda^3}$&$\frac{1}{4\lambda}+\frac{1}{B\lambda^2}+\frac{1}{2B^2\lambda^3}$\\
$t^3$&$ \frac12-\frac{9}{4A\lambda}+\frac{3}{2A^3\lambda^3}$&$~ \frac12+\frac{9}{4B\lambda}-\frac{3}{2B^3\lambda^3}~$&$ \frac{\lambda}{4}-\frac{3}{\lambda t^2}$\rule{0mm}{5mm}
&$~\frac{1}{4\lambda}-\frac{3}{2A\lambda^2}+\frac{3}{2A^2\lambda^3}~$&$\frac{1}{4\lambda}+\frac{3}{2B\lambda^2}+\frac{3}{2B^2\lambda^3}$\\
$t^4$&$ \frac12-\frac{3}{A\lambda}+\frac{6}{A^3\lambda^3}$&$ \frac12+\frac{3}{B\lambda}-\frac{6}{B^3\lambda^3}$&$ \frac{\lambda}{4}-\frac{6}{\lambda t^2}$\rule{0mm}{5mm}
&$\frac{1}{4\lambda}-\frac{2}{A\lambda^2}+\frac{3}{A^2\lambda^3}$&$\frac{1}{4\lambda}+\frac{2}{B\lambda^2}+\frac{3}{B^2\lambda^3}$\\
$e^t$&$ \frac12-\frac{3}{4\lambda}+\frac{1}{4\lambda^3}$&$ \frac12+\frac{3}{4\lambda}-\frac{1}{4\lambda^3}$&$ \frac{\lambda}{4}-\frac{1}{2\lambda }$\rule{0mm}{5mm}
&$\frac{1}{4\lambda}-\frac{1}{2\lambda^2}+\frac{1}{4\lambda^3}$&$\frac{1}{4\lambda}+\frac{1}{2\lambda^2}+\frac{1}{4\lambda^3}$\\
\\
\hline
\end{tabular}
\label{tab:optdes-ar2}
\end{table}

\subsection{Practical implementation}

Suppose that the  $N$ equidistant points defined in \eqref{points} are the potential observation points.
Let $K+2$ be the number of observations actually taken in the experiment and that we want to
construct  a discrete design, which  can  be implemented in practice.
Suppose that $K$ is small and $N$ is large, then  efficient designs and corresponding estimators for the model \eqref{one}
can be  derived
from the continuous approximations,  which have been developed in the previous sections.

In \cite{DetPZ2016} a procedure with a  good finite sample performance is proposed.
It consists of a slight modification of the SLSE given in \eqref{eq:slse}
and a discretization of the density $p(t)$ defined in \eqref{eq:ar1-pt} for AR(1) errors
and \eqref{eq:ar2-ptgen} for AR(2) errors. To be precise  consider a continuous SLSE with weights at the points $A$ and $B$
(the end-points of the interval $[A,B]$), which correspond to the masses $P_A$ and $P_B$ and,
for the AR(2) errors, $Q_A$ and $Q_B$ as well.
We thus only need to approximate the continuous part of the design, which has a density on $(A,B)$,
by a $K$-point design with equal masses.

We assume that the density $p(\cdot )$ is not identically zero on the interval $(A,B)$.
Define
$\varphi(t) = \kappa|p(t)|$ for $ t \in (A,B)$
 and
choose the constant $\kappa$ such that $\int_A^B \varphi(t)dt =1$,
that is,
\bea
 \kappa=\frac{1}{\int_A^B |p(t)| dt}.
\eea
 Denote by
$F(t)= \int_A^t \varphi(s) ds$  the  corresponding cumulative distribution function.
As $K$-point design we use a $K$-point approximation to the measure with density $\varphi(t)$, that is
$\hat{\xi}_{K}=\{ t_{1,K}, \ldots,t_{K,K}; 1/K, \ldots, 1/K\}  $,
where $t_{i,K}= R(F^{-1}(i/(K+1)))$     $i=1,2, \ldots, K$.
Here $R(t)$ is the operator of rounding a number $t$ towards the set of   points defined by \eqref{points},
 that is $R(F(i/(K+1))= t_{i,K} := A+ (i^*-1)\Delta$, where
 $$
 |F(i/(K+1)) - A + (i^*-1) \Delta | = \min \{ | F(i/(K+1))- A+(j-1)\Delta |\, ; ~~ j=1,\ldots,N \}.
 $$

If $p(t)=0 $ on a sub-interval of $[A,B]$  and $ F^{-1}(i/(K+1))$ is not uniquely defined then  we choose the smallest element
from the set $ R(F^{-1}(i/(K+1))$ as $t_{i,K}$.
Also we define $s_{i,K}= {\rm sign} (p(t_{i,K}))$ and obtain
from the representation of the continuous SLSE for AR(1) errors  in  Proposition~\ref{prop1}
a reasonable estimator with corresponding design. To be precise,
  $y_1,\ldots , y_{K+2}$  should be observed at experimental conditions
$A, t_{1,K}, t_{2,K}, \ldots,$ $t_{K,K}, B$, respectively,  and   the parameter $\theta$
has to be estimated by the {\it modified} SLSE
\bea
 \hat\theta_{K+2}\!=\!D_{K+2} \Big(\!P_Af(A)y_{A}\!+\!P_Bf(B)y_{B}
 \!+\!\frac{B-A}{\kappa K}\sum_{i=1}^K \!\!s_{i,K}  f(t_{i,K})y_i
 \Big),
\eea
where
\bea
 D_{K+2}\!=\! \Big(P_A f^2(A)\!+\!P_B f^2(B)
  \!+\!\frac{B-A}{\kappa K}\sum_{i=1}^K s_{i,K}  f^2(t_{i,K})\Big)^{-1}.
\eea
It follows from the discussion of the previous paragraph that
$ \mathrm{Var} (\hat \theta_{K+2}) \approx D^*$,
where $D^*$ is defined in \eqref{eq:modif-slse-ar1}.
Similarly, the modified  SLSE for AR(2) errors  is
defined by
{\small
\be
 \hat\theta_{K+2}\!=\!D_{K+2}\Big(\!\!\!\!\!\!&&Q_Bf(B)y'(B)\!-\!Q_Af(A)y'(A)\!\nonumber \\
 &&+\!P_Af(A)y_{A}\!+\!P_Bf(B)y_B
 \!+\!\frac{B-A}{\kappa K}\sum_{i=1}^K s_{i,K} f(t_{i,K})y(t_{i,K}) \Big)
\label{eq:modif-slse-K-ar2}
\ee}
where
\bea
 D_{K+2}\!=\! \Big(Q_Bf(B)f'(B)\!-\!Q_Af(A)f'(A)\!+\!P_Af^2(A)\!+\!P_Bf^2(B)
  \!+\!\frac{B\!-\!A}{\kappa K}\sum_{i=1}^K s_{i,K}  f^2(t_{i,K})\Big)^{-1}.
\eea
In \eqref{eq:modif-slse-K-ar2}  the expressions are the derivatives   $y'(A)$ and $y'(B)$  of the continuous approximation
$\{y(t) \}_{t\in [A,B]}$, which are usually not available in practice.
Therefore, we recommend to make two additional observations at the points
$A+\Delta$ and $B-\Delta$ and  to replace
the derivatives by their approximations  $(y_{A+\Delta}-y_{A})/\Delta$ and $(y_{B}-y_{B-\Delta})/\Delta$.
Thus, we replace the estimator \eqref{eq:modif-slse-K-ar2}
by   the weighted least squares estimator (WLSE)
\be \label{wlsemod}
 \tilde \theta_{K+4}  = ({X}^T W {X})^{-1}{X}^T W {Y},
\ee
where $Y=(y_A,y_{A+\Delta},y_{t_{1,K}}, \ldots,y_{t_{K,K}},y_{B-\Delta},y_{B})^T$
and the matrix $W$ is defined  by
\be \label{wlse}
 W=\mathrm{diag}\Big\{
 \frac{P_A}2+\frac{Q_A}\Delta,\frac{P_A}2-\frac{Q_A}\Delta,s_{1,K}\frac{B\!-\!A}{\kappa K}, \ldots,
 s_{K,K}\frac{B\!-\!A}{\kappa K},\frac{P_B}2-\frac{Q_B}\Delta,\frac{P_B}2+\frac{Q_B}\Delta
 \Big\}.
\ee
Note that the variance of $\tilde \theta_{K+4}$ is given by
\bea
 \mathrm{Var}(\tilde \theta_{K+4})=({X}^TW{X})^{-1} (X^TW\Sigma^{}WX) ({X}^TW{X})^{-1}.
\eea

\subsection{Practical performance}

Consider the regression model \eqref{eq:model-onepar} with $f(t)=1$, $[A,B]=[0,1]$
and AR(2) errors. Suppose that $N=101$ so that $t_i=i/100$, $i=0,1,\ldots,N$, are potential observation points.
We also assume that the autocorrelation function $r_k$ is of the form \eqref{eq:K-elin} with
$\lambda=1$.
We investigate the  design $\xi_{K+2}$ with $(K+2)$ points $0, t_{1,K}, t_{2,K}, \ldots,$ $t_{K,K}, 1$ and
  the  design $\xi_{K+4}$  with $(K+4)$  points $ 0,0.01, t_{1,K}, t_{2,K}, \ldots,$ $t_{K,K}, 0.99,1$.
The points $ t_{1,K}, t_{2,K}, \ldots,$ $t_{K,K} $ are shown in  the second column of Table \ref{tab:example0}.
 In this table we  also display  the  variances
of the WLSE   $\tilde \theta_{K+4}$,  defined by \eqref{wlse}, the LSE $\hat \theta_{\LSE,K+2}$ based on the design
$\xi_{K+2} $  and the BLUE $\hat \theta_{\BLUE,K+2}$ and $\hat \theta_{\BLUE,K+4}$   for the  designs
$\xi_{K+2} $ and $\xi_{K+4} $, respectively.
Let $\hat \theta_{\BLUE}$ denote the BLUE based on 101 observations at the points $\{ \frac {i}{100} | \ i=0,\ldots,100 \}$, then we
observe  that $0.80158449 = \mathrm{Var}(\hat \theta_{ \mathrm{\BLUE}})\approx D^*=0.8$
that is in agreement with Theorem \ref{th:approx-ar2}.
We also observe that $\mathrm{Var}(\hat \theta_{\BLUE,K+4})\cong \mathrm{Var}(\hat \theta_{ \mathrm{BLUE}})$
and $\mathrm{Var}(\hat \theta_{\BLUE,K+2})\not\cong \mathrm{Var}(\hat \theta_{ \mathrm{\BLUE}})$
showing the importance of taking one additional   observation at  each boundary point $A$ and $B$.
Note that the proposed estimator $\tilde \theta_{K+4}$ defined in \eqref{wlse}  is nearly as accurate as the
BLUE $\hat \theta_{\BLUE,K+4}$ at the same points
and that the LSE $\hat \theta_{\LSE,K+2}$ is about $ 10-15\% $ worse than the BLUE.

\begin{table}[h]
\caption{\it The variances of the LSE, the WLSE defined by \eqref{wlse}
and the BLUE for designs with $K+2$ and $K+4$ points. $f(t)=1$, $[A,B]=[0,1]$, $N=101$,
the autocovariance structure  is given by \eqref{eq:K-elin} with $\lambda=1$, which yields
 $D^*=0.80000$ and
$\mathrm{Var}(\hat \theta_{\small \rm BLUE})=0.80158449$.
}
\begin{center}
\begin{tabular}{cccccc}
\hline
 $K$&$t_{1,K}, \ldots,t_{K,K}$&$\mathrm{Var}(\hat \theta_{\LSE,K+2})$& $\mathrm{Var}(\tilde \theta_{K+4})$
 & $\mathrm{Var}(\hat \theta_{\BLUE,K+2})$ & $\mathrm{Var}(\hat \theta_{\BLUE,K+4})$\\
\hline
 $2$&0.33, 0.67        &0.914 &0.80170&0.82663 &0.80158714\\
 $3$&0.25, 0.5, 0.75   &0.921 &0.80165&0.82022 &0.80158533\\
 $4$&0.2, 0.4, 0.6, 0.8&0.925 &0.80162&0.81681 &0.80158484\\
 5&0.17, 0.33, 0.5, 0.67, 0.83&0.928 &0.80161&0.81443 &0.80158466\\
\hline
\end{tabular}
\end{center}
\label{tab:example0}
\end{table}

\vspace{2mm}
As a second example,
consider the regression model \eqref{eq:model-onepar} with $f(t)=t^2$, $[A,B]=[0.1,1.1]$
and AR(2) errors. Suppose that $N=101$ so that $t_i=0.1+i/100$, $i=0,1,\ldots,N$, are potential observation points.
We also assume that the autocorrelation function $r_j$ is of the form \eqref{eq:K-elin} with
$\lambda=2$.
We investigate the design $\xi_{K+2}$ with $(K+2)$ points $0.1, t_{1,K}, t_{2,K}, \ldots,$ $t_{K,K}, 1.1$ and the
design $\xi_{K+4}$  with $(K+4)$ points $0.1,0.11,$ $t_{1,K}, t_{2,K}, \ldots, t_{K,K}, 1.09,1.1$.
The non-trivial points are  shown  in  the second column of  Table \ref{tab:example2}. In the other
columns we display   the  variances of the different estimators introduced in the previous paragraph.
We observe again  that $0.37055791= \mathrm{Var}(\hat \theta_{ \BLUE})\approx D^* =0.36543 $
that is in line with Theorem \ref{th:approx-ar2}.
Note  also  that $\mathrm{Var}(\hat \theta_{\BLUE,K+4})\cong \mathrm{Var}(\hat \theta_{ \BLUE})$
and the estimator $\hat \theta_{\BLUE,K+2}$ without the two additional observations at the boundary is not efficient.
Again the proposed estimator $\tilde \theta_{K+4}$ is nearly as accurate as the BLUE at the same points
but the LSE $\hat \theta_{\LSE,K+2}$ is dramatically worse than the BLUE.

\begin{table}[!hhh]
\caption{\it The variances of the LSE, the WLSE
and the BLUE for designs with $K+2$ and $K+4$ points.  $f(t)=t^2$, $[A,B]=[0.1,1.1]$, $N=101$
and the   autocovariance is given by  \eqref{eq:K-elin} with $\lambda=2$, which
yields $D^*=60000/164189\cong 0.36543$ and
$\mathrm{Var}(\hat \theta_{ \mathrm{BLUE}})=0.37055791$.
}
\begin{center}
\begin{tabular}{cccccc}
\hline
 $K$&$t_{1,K}, \ldots,t_{K,K}$&$\mathrm{Var}(\hat \theta_{\LSE,K+2})$& $\mathrm{Var}(\tilde \theta_{K+4})$
 & $\mathrm{Var}(\hat \theta_{\BLUE,K+2})$ & $\mathrm{Var}(\hat \theta_{\BLUE,K+4})$\\
\hline
 $2$&0.14, 0.22                & 0.723 &0.40218 &0.53175 &0.37079053 \\
 $3$&0.12, 0.17, 0.27          & 0.751 &0.40204 &0.52509 &0.37072082\\
 $4$&0.12, 0.15, 0.20, 0.30    & 0.783 &0.40176 &0.52089 &0.37068565\\
 5&0.12, 0.14, 0.17, 0.22, 0.33& 0.818 &0.40139 &0.51689 &0.37065785\\
\hline
\end{tabular}
\end{center}
\label{tab:example2}
\end{table}

\medskip
\medskip

{\bf Acknowledgements.}
This work has been supported in part by the
Collaborative Research Center "Statistical modeling of nonlinear
dynamic processes" (SFB 823, Project C2) of the German Research Foundation (DFG) and the
   National Institute Of General Medical Sciences of the National
Institutes of Health under Award Number R01GM107639. The content is solely the responsibility of the authors and does not necessarily
 represent the official views of the National
Institutes of Health.
The work of Andrey Pepelyshev was partly supported by
the project "Actual problems of design and analysis for regression models" (6.38.435.2015)
of St. Petersburg State University.
The authors would like to thank Martina
Stein, who typed parts of this manuscript with considerable
technical expertise.
% We are very grateful to a referee, the associate editor and the editor for their
% constructive comments, which led to substantial improvement of  an earlier version of this manuscript.
\bibliography{opt_designs}

\end{document}